\documentclass{article}

\newcommand{\rom}{\rm}

\newtheorem{theorem}{Theorem}[section]
\newtheorem{lemma}[theorem]{Lemma}

\newtheorem{remark}[theorem]{Remark}
\newtheorem{example}[theorem]{Example}

\title
{The Calibration Method\\ for Free Discontinuity Problems}
\author{Gianni Dal Maso
\thanks {SISSA, Via Beirut 4, 34014 Trieste, Italy, 
e-mail address: {\tt dalmaso@sissa.it}}
}
\date{}

\begin{document}

\maketitle
\begin{abstract}
The calibration method is used to identify some
minimizers of the Mumford-Shah functional. 
The method is then extended to more general free discontinuity 
problems.
\end{abstract}

\section{Introduction}
In \cite{DG} De Giorgi introduced the name
{\it free discontinuity problems\/} to denote a wide class of
minimum problems for functionals of the form
\begin{equation}\label{fg}
F(u):=\int_{\Omega\setminus S_u} f(x,u(x),\nabla u(x))dx \,
+ \int_{S_{u}} \psi(x,u^{+}\!(x),u^{-}\!(x),\nu_{u}(x))
d{\cal H}^{n-1}\!,
\end{equation}
where $\Omega$ is a given bounded domain in 
${\mathbf R}^n$ with Lipschitz 
boundary,
$f\colon\Omega\times{\mathbf R}
\times{\mathbf R}^{n}\to[0,+\infty]$ and 
$\psi\colon\Omega\times{\mathbf R} 
\times{\mathbf R}\times{\mathbf S}^{n-1}\to[0,+\infty]$ 
are given Borel functions, ${\mathbf S}^{n-1}=
\{{v\in{\mathbf R}^{n}}:{|v|=1}\}$,
${\cal H}^{n-1}$ is the $({n-1})$-dimensional Hausdorff measure,
and the unknown function $u\colon \Omega\to{\mathbf R}$ is assumed
to be regular out of a (partially regular) singular set $S_u$ of dimension 
${n-1}$, with unit normal $\nu_{u}$, on which $u$
admits unilateral traces $u^{+}$ and $u^{-}$.  The main feature of 
these problems is that the shape and location of 
the discontinuity set $S_{u}$ are not prescribed.  
Thus minimizing $F$ means optimizing both the 
function  $u$ and the singular set $S_{u}$, 
which is indeed often regarded  
as an independent unknown.
\par
These problems have an 
increasing importance in many branches of applied analysis, such as 
image processing (Mumford-Shah functional for image segmentation) 
and  fracture mechanics (Griffith's criterion and 
Barenblatt cohesive zone model).
\par
The {\it Mumford-Shah functional\/} was introduced in \cite{MS} 
in the context of a variational approach to image segmentation problems
(for which we refer to \cite{MoS}).
It can be written, in dimension $n$, as
\begin{equation}\label{1}
F^{\alpha,\beta}_{g}(u):=
\int_{\Omega \setminus S_u} \hskip-7pt | \nabla u(x)|^2dx 
+ \alpha  {\cal H}^{n-1}(S_u) 
+ \beta \int_{\Omega \setminus S_u} \hskip-7pt |u(x)-g(x)|^2 dx \ ,
\end{equation}
where $g$
is a given function in $L^{\infty}(\Omega)$ 
(interpreted as the grey level of the image to be 
analysed), and $\alpha>0$ and $\beta\ge 0$ are constants.
When $n=2$ (the only case considered in image 
processing), the singular set $S_{u}$ of a minimizer $u$ of
$F^{\alpha,\beta}_{g}$ is interpreted as the 
set of the most relevant segmentation lines of the image.
\par
Using different classes of infinitesimal variations, 
one can show that every minimizer must satisfy certain
equilibrium conditions, which could be globally called
{\it Euler-Lagrange equations\/} for $F^{\alpha,\beta}_{g}$.
For instance, $u$ must 
satisfy the equation $\Delta u=\beta({u-g})$ 
on $\Omega\setminus S_u$, with Neumann boundary conditions
on $S_u\cup\partial\Omega$.
Moreover, there is a link between the mean curvature of $S_u$
(where defined)
and the traces of $u$ and $\nabla u$ on the two sides of $S_u$;
for instance, when $\beta=0$, the mean curvature of $S_u$ must be 
equal to the difference of the squares of the norms of the traces of $\nabla u$.
Additional conditions have been derived for the 
two-dimensional case. We refer the reader to \cite{MS} and 
\cite{AFP}
for a precise description of these equilibrium conditions.
\par
However, since $F^{\alpha,\beta}_{g}$ is not convex,
all conditions which can be derived by infinitesimal
variations are necessary for minimality, but never sufficient. 
The purpose of this note is precisely to present a {\it sufficient 
condition\/} for minimality 
(Theorem~\ref{3.2} for $F^{\alpha,\beta}_{g}$ and Theorem~\ref{03.1}
for $F$), and give a few 
applications (Examples \ref{ex1}--\ref{ex8}). Detailed proofs and further 
results will be given in the forthcoming paper \cite{ABD}.

\section{Notation and preliminaries}
For a complete mathematical treatment of the minimum problems 
for the functional $F$ considered in (\ref{fg}),
we use the space $SBV(\Omega)$ of {\it special functions of bounded 
variation\/}, introduced by De Giorgi and Ambrosio in \cite{DG-Amb}.
A self-con\-tained presentation of this space can be found in the recent 
book \cite{AFP}, which contains also the complete proof of the 
existence of a minimizer $u$ of $F^{\alpha,\beta}_{g}$, and of
the partial regularity of the corresponding singular set $S_{u}$ (the 
regularity of $u$ on $\Omega\setminus S_{u}$ follows from the standard 
theory of elliptic equations).
\par
We recall that for every $u\in SBV(\Omega)$ the 
{\it approximate upper and 
lower limits\/} $u^{+}\!(x)$ and $u^{-}\!(x)$ 
at a point $x\in\Omega$ are defined by
$$
u^{\pm}\!(x):=\pm \inf\,\{t\in{\mathbf R}:\lim_{\rho\to 0+} \,
\rho^{-n} {\cal L}^{n}(\{\pm u>t\}\cap 
B_{\rho}(x))=0\} \ ,
$$
where $B_{\rho}(x)$ is the open ball with centre $x$ and radius $\rho$.
The {\it singular set\/} (or {\it jump set\/}) of  $u$ is defined by
$S_u:=\{{x\in \Omega}: u^{-}\!(x)<u^{+}\!(x)\}$.
It is known that $S_u$ is countably 
$({\cal H}^{n-1},n-1)$-rectifiable and that there exists a Borel measurable 
function $\nu_{u}\colon S_u\to {\mathbf S}^{n-1}$ such that  for 
${\cal H}^{n-1}$-a.e.\ $x\in S_u$ we have
\begin{equation}\label{am1}
\lim_{\rho\to0+}\,{1\over\rho^{n}}
\int_{B^{\pm}_{\rho}(x)} |u(y)-u^{\pm}\!(x)|\,dy = 0\ ,
\end{equation}
where $B^{\pm}_{\rho}(x):=\{y\in B_{\rho}(x): 
\pm(y-x)\cdot\nu_{u}(x)>0\}$ and $\cdot$ denotes the scalar product 
in ${\mathbf R}^{n}$ (see \cite[Theorem~4.5.9]{Fed}). 
Condition (\ref{am1}) says that $\nu_{u}(x)$ 
points from the side of $S_u$ corresponding to $u^{-}\!(x)$ to the side 
corresponding to $u^{+}\!(x)$.
\par
The gradient $Du$ of $u$ is a measure that can be decomposed 
as the sum of two measures
$
Du=D^{a}u + D^{s}u
$,
where $D^{a}u$ is absolutely continuous and
$D^{s}u$ is singular with respect to the Lebesgue measure ${\cal L}^{n}$.
The density of $D^{a}u$ with respect to ${\cal L}^{n}$ is denoted by
$\nabla u$. Since $u\in SBV(\Omega)$, for every Borel set 
$B$ in $\Omega$ we have
$$
(Du)(B) = \int_{B} \nabla u(x)\,dx +
\int_{B\cap S_u} \hskip-5pt (u^{+}\!(x)-u^{-}\!(x))\,
\nu_{u}(x)\,d{\cal H}^{n-1} \,.
$$
\par
The {\it graph\/} of $u$ is defined as
$$
{\mit\Gamma}_u := \{(x,t)\in \Omega\times{\mathbf R} : u^{-}\!(x)\le t 
\le u^{+}\!(x)\} \ .
$$
The characteristic function of the subgraph 
$\{{(x,t)\in \Omega\times{\mathbf R} }: {t \le u(x)}\}$ is denoted by $1_{u}$. 
It is defined  by $1_u(x,t):=1$ if $t\le u(x)$, and $1_u(x,t):=0$ if $t>u(x)$.
It belongs to 
$SBV({\Omega \times{\mathbf R}})$ and its gradient $D1_u$ is a measure
concentrated on ${\mit\Gamma}_u$.

\section{The main results}
We fix an open subset $U$ of
${\Omega\times{\mathbf R}}$ of the form
\begin{equation}\label{U}
U:=\{(x,t)\in {\Omega\times{\mathbf R}}: \tau_{1}(x)<t<\tau_{2}(x)\}\ ,
\end{equation}
where $\tau_{1}$ and $\tau_{2}$ are two continuous functions on 
$\overline\Omega$
such that
$-\infty\le \tau_{1}(x) \le \tau_{2}(x)\le+\infty$ for every 
$x\in\overline \Omega$.
\par
Let $F$ be the functional introduced in (\ref{fg}).
We say that a function $u\in SBV(\Omega)$, 
with graph ${\cal H}^{n}$-contained in $U$ (i.e., 
${\cal H}^{n}({{\mit\Gamma}_u
\setminus U})=0$),
is a {\it Dirichlet $U$-minimizer\/} of $F$, 
if $F(u)\le F(v)$ for every 
$v\in SBV(\Omega)$ with the same trace as $u$ on $\partial\Omega$ and 
with graph ${\cal H}^{n}$-contained in $U$.  
If the inequality $F(u)\le F(v)$ 
holds for every 
$v\in SBV(\Omega)$ with graph ${\cal H}^{n}$-contained in $U$, 
we say that $u$ 
is a {\it $U$-minimizer\/} of $F$. We omit 
$U$ when $U={\Omega\times {\mathbf R}}$.
\par
The symbol $\phi$ will always denote
a bounded Borel measurable vectorfield defined on 
$U$ with values in ${\mathbf R}^{n+1}={\mathbf R}^{n}\times{\mathbf R}$,
with components
$\phi^x\in{\mathbf R}^n$ and $\phi^t\in{\mathbf R}$. The divergence of $\phi$ is then
${\rm div}\phi(x,t)={\rm div}_{x}\phi^{x}(x,t) +\partial_{t}\phi^t(x,t)$.
\par
We begin with a theorem concerning the functional 
$F^{\alpha,\beta}_{g}$ introduced in (\ref{1}).

\begin{theorem}\label{3.2}
Let $u\in SBV(\Omega)$ with graph ${\cal H}^{n}$-contained in $U$.
Assume that there exists a bounded vectorfield
$\phi$ of class $C^1$ on $U$ with the following properties:
\begin{description}
\smallskip
\item[$\rm(a1)$]
\quad${1\over 4} |\phi^x(x,t)|^2 \le \phi^t(x,t)+ \beta|t-g(x)|^2$
\par\vspace{-5pt}
\item[\hphantom{$\rm(a1)$}]
\quad\qquad for ${\cal L}^{n}$-a.e.\ $x\in\Omega$ and for
every $\tau_{1}(x)<t <\tau_{2}(x)$;
\smallskip
\item[$\rm(a2)$]
\quad$\phi^x(x,u(x))=2\,\nabla u(x)\ $ and
$\ \phi^t(x,u(x))=|\nabla u(x)|^2-\beta|u(x)-g(x)|^2$ 
\par\vspace{-5pt}
\item[\hphantom{$\rm(a2)$}]
\quad\qquad  for ${\cal L}^{n}$-a.e.\ $x\in\Omega$;
\smallskip
\item[$\rm(b1)$]
\quad$\displaystyle\Big| \int_{t_1}^{t_2} \phi^x(x,t) \, dt 
\Big|\le\alpha$
\par\vspace{-5pt}
\item[\hphantom{$\rm(b1)$}]
\quad\qquad for ${\cal H}^{n-1}$-a.e.\ $x\in\Omega$  and for
every $\tau_{1}(x)<t_{1}<t_{2} <\tau_{2}(x)$;
\smallskip
\item[$\rm(b2)$]
\quad$\displaystyle\int_{u^-(x)}^{u^+(x)} \phi^x(x,t) \, dt
=\alpha\,\nu_u(x)$ 
\quad for ${\cal H}^{n-1}$-a.e.\ $x\in S_u$;
\smallskip
\item[$\rm(c1)$]
\quad$\displaystyle {\rm div}\phi(x,t)=0$ \quad for every $(x,t)\in U$.
\smallskip
\end{description}
\noindent
Then $u$ is a Dirichlet $U$-minimizer of $F^{\alpha,\beta}_{g}$. 
If, in addition, $\phi^{x}(x,t)$ satisfies 
the boundary condition
\begin{description}
\smallskip
\item[$\rm(c2)$]
\quad$\displaystyle\lim_{(y,s)\to (x,t)} 
\phi^x(y,s)\cdot \nu(x)=0$
\par\vspace{-5pt}
\item[\hphantom{$\rm(c2)$}]
\quad\qquad for ${\cal H}^{n-1}$-a.e.\ $x\in\partial\Omega$ 
and for ${\cal L}^{1}$-a.e.\ $t\in 
[\tau_{1}(x),\tau_{2}(x)]$,
\smallskip
\end{description}
\noindent
where $\nu(x)$ is the outer unit normal 
to $\partial\Omega$,
then $u$ is a $U$-minimizer of $F^{\alpha,\beta}_{g}$.
\end{theorem}

A vectorfield $\phi$ which satisfies conditions $\rm(a1)$--$\rm(c1)$
of Theorem~\ref{3.2} is called a {\it calibration\/} for the 
functional $F^{\alpha,\beta}_{g}$ on $U$. If $\phi$ satisfies also
$\rm(c2)$, it is called a  {\it Neumann calibration\/}.
Theorem \ref{3.2} is an immediate consequence of 
the following lemmas.

\begin{lemma}\label{3.7}
Let $\phi$ be
a vectorfield which satisfies conditions $\rm(a1)$ and $\rm(b1)$
of Theorem~\ref{3.2}.
Then for every $u\in SBV(\Omega)$ 
with graph ${\cal H}^{n}$-contained in $U$
we have
\begin{equation}\label{(3.7}
F^{\alpha,\beta}_{g}(u)\ge \int_{U}  \phi\cdot d(D1_u) \ .
\end{equation}
Moreover, equality holds in (\ref{(3.7}) for a given 
$u$ if and only if conditions $\rm(a2)$ and $\rm(b2)$ of 
Theorem \ref{3.2} are satisfied.
\end{lemma}

The next lemma is a consequence of the divergence theorem.

\begin{lemma}\label{ae1}
Suppose that  $\phi$ is of class $C^{1}$ and that 
${\rm div}\phi=0$ on $U$.
Then
\begin{equation}\label{ab44}
\int_{U} \phi\cdot d(D1_u) = \int_{U} \phi\cdot d(D1_v) 
\end{equation}
for every pair of functions $u$, $v$ in $BV(\Omega)$ 
with the same trace on $\partial \Omega$ and
with graphs ${\cal H}^{n}$-contained in $U$. If, in addition, $\phi$
satisfies
condition $\rm(c2)$ of Theorem~\ref{3.2}, then
(\ref{ab44}) holds for every pair of functions $u$, $v$ in $BV(\Omega)$
with graphs ${\cal H}^{n}$-contained in $U$.
\end{lemma}

As a matter of fact, the method of calibrations can
be easily adapted to the functional $F$ defined in (\ref{fg}).

\begin{theorem}\label{03.1}
Let $u\in SBV(\Omega)$ with graph ${\cal H}^{n}$-contained in $U$.
Assume that there exists a bounded vectorfield
$\phi$ of class $C^1$ on $U$ with the following properties:
\begin{description}
\smallskip
\item[$\rm(a1)$]
\quad$\phi^x(x,t)\cdot v\le \phi^t(x,t) + f(x,t,v)$
\qquad for ${\cal L}^{n}$-a.e.\ $x\in\Omega$,
\par\vspace{-5pt}
\item[\hphantom{$\rm(a1)$}]
\quad\qquad for every
$\tau_{1}(x)<t <\tau_{2}(x)$,
and for every $v\in{\mathbf R}^{n}$;
\smallskip
\item[$\rm(a2)$]
\quad$\phi^x(x,u(x))\cdot \nabla u(x) = 
\phi^t(x,u(x)) + f(x,u(x),\nabla u(x))$
\par\vspace{-5pt}
\item[\hphantom{$\rm(a2)$}]
\quad\qquad for ${\cal L}^{n}$-a.e.\ $x\in\Omega$;
\smallskip
\item[$\rm(b1)$]
\quad$\displaystyle \nu\cdot\int_{t_1}^{t_2} \phi^x(x,t) \, dt 
\le\psi(x,t_1,t_2,\nu)$\qquad for ${\cal H}^{n-1}$-a.e.\ $x\in\Omega$, 
\par\vspace{-5pt}
\item[\hphantom{$\rm(b1)$}]
\quad\qquad 
for every $\tau_{1}(x)<t_{1}<t_{2} <\tau_{2}(x)$, and for every
$\nu\in {\mathbf S}^{n-1}$;
\smallskip
\item[$\rm(b2)$]
\quad$\displaystyle \nu_u(x) \cdot \int_{u^-\!(x)}^{u^+\!(x)} 
\phi^x(x,t) \, dt
 = \psi ( x,u^-(x),u^+(x),\nu_u(x))$ 
\par\vspace{-5pt}
\item[\hphantom{$\rm(b2)$}]
\quad\qquad for ${\cal H}^{n-1}$-a.e.\ $x\in S_u$;
\smallskip
\item[$\rm(c1)$]
\quad$\displaystyle {\rm div}\phi(x,t)=0$ \quad for every $(x,t)\in U$.
\smallskip
\end{description}
\noindent
Then $u$ is a Dirichlet $U$-minimizer of $F$. 
If $\phi^{x}(x,t)$ satisfies also
the boundary condition $\rm(c2)$ of Theorem~\ref{3.2},
then $u$ is a $U$-minimizer of $F$.
\end{theorem}

\begin{remark} \rom
We note that in Theorem~\ref{03.1} there is no regularity or 
convexity hypothesis on $f$ or $\psi$.
If $f^*(x,t,v^{*})$ is the the convex 
conjugate of $f(x,t,v)$ with respect to $v$, condition $\rm(a1)$ is 
equivalent to
\begin{description}
\smallskip
\item[$\rm(a1')$]
\quad$f^*(x,t,\phi^x(x,t)) \le \phi^t(x,t)$
\par\vspace{-5pt}
\item[\hphantom{$\rm(a1')$}]
\quad\qquad for ${\cal L}^{n}$-a.e.\ $x\in\Omega$ and for every
$\tau_{1}(x)<t <\tau_{2}(x)$.
\smallskip
\end{description}
\noindent
If this condition is satisfied, and
$f(x,t,v)$ is convex and differentiable with respect to $v$, then
condition $\rm(a2)$ is 
equivalent to
\begin{description}
\smallskip
\item[$\rm(a2')$]
\quad$\cases{\phi^x(x,u(x))=
\partial_{v} f(x,u(x),\nabla u(x))&\cr
\vphantom{\vrule 
width 0pt height 13pt depth 0pt}
\phi^t(x,u(x))=f^*(x,u(x),\phi^x(x,u(x)))&\cr}$
\ for ${\cal L}^{n}$-a.e.\ $x\in\Omega$.
\end{description}
\noindent
\end{remark}

\begin{remark}\label{rem}  \rom
In Theorems~\ref{3.2} and~\ref{03.1}
the hypothesis that $\phi$ is of class $C^{1}$ 
is too strong 
for many applications. It is used only in Lemma~\ref{ae1} and it can 
be relaxed
in several ways (see \cite{ABD} for details). 
For instance, one may consider piecewise $C^1$ vectorfields, 
which may be discontinuous along sufficiently regular
interfaces. In this case the divergence-free condition  $\rm(c1)$
must be understood in the distributional sense, i.e., 
the pointwise divergence  vanishes (where defined) and 
the normal component of $\phi$ is continuous across 
the discontinuity surfaces.
\end{remark}

\section{Some examples}
The following examples show that the calibration method is very 
flexible, and can be used to prove the minimality of a given function 
$u$ in many different situations. In the first examples
we will consider only the ``homogeneous'' functional 
$F^{\alpha}:=F^{\alpha,0}_{g}$, 
in which the lower order term $\beta\int_{\Omega}|{u-g}|^{2} dx$
vanishes.

\begin{example}[Affine function in one dimension]\label{ex1} \rom
Let $n:=1$, $\Omega:={]0,a[}\,$, and $u(x):=\lambda x$, 
with $\lambda>0$. 
It is easy to see that
$u$ is a Dirichlet minimizer of $F^{\alpha}$ if and only if 
$a\lambda^{2}\le \alpha$. In this case a calibration is given by the 
piecewise constant function
\begin{equation}\label{6}
\phi(x,t):=\cases{(2\lambda, \lambda^{2}),&if
$\frac{\lambda}2x\le t \le \frac{\lambda}2(x+a)$,
\cr\vphantom{\vrule 
width 0pt height 13pt depth 0pt}
(0,0),&otherwise.
\cr}
\end{equation}
Another calibration is given by
\begin{equation}\label{6b}
\phi(x,t):=\cases{\big(2\frac{t}{x}, (\frac{t}{x})^{2}\big),&if
$0\le t\le\lambda x$,
\cr\vphantom{\vrule width 0pt height 13pt depth 0pt}
\big(2\frac{\lambda a - t}{a- x}, (\frac{\lambda a - t}{a- x})^{2}\big),&if
$\lambda x\le t\le\lambda a$,
\cr\vphantom{\vrule width 0pt height 13pt depth 0pt}
(0,0),&otherwise.
\cr}
\end{equation}
\par
If $a\lambda^{2}> \alpha$, then the function $u(x):=\lambda x$ is not 
a Dirichlet minimizer of $F^{\alpha}$, but it is still a Dirichlet 
$U$-minimizer with
$$
\textstyle U:=\{(x,t)\in{]0,a[}\times{\mathbf R}: 
\lambda x - \frac{\alpha}{4\lambda} < t 
<\lambda x +\frac{\alpha}{4\lambda}\}\,.
$$
A calibration on $U$ is given by $\phi(x,t):=(2\lambda, \lambda^{2})$.
\end{example}

\begin{example}[Jump in one dimension]\label{ex2} \rom
Let $n:=1$, $\Omega:={]0,a[}\,$, $u(x):=0$ for $0<x<c$, and $u(x):=h$ for
$c<x<a$, with $0<c<a$ and $h>0$. It is easy to see that
$u$ is a Dirichlet minimizer of $F^{\alpha}$ if and only if 
$a\alpha \le h^{2}$. In this case two different calibrations are given by 
(\ref{6}) and (\ref{6b})  with $\lambda=\sqrt\alpha/\sqrt{a}$.
\par
Suppose now that $a\alpha > h^{2}$. Let $\varepsilon>0$ be a 
constant such that
$2\varepsilon + \sqrt{2\alpha\varepsilon}\le h$, let
$$
\tau_{1}(x)=\cases{-\varepsilon,&if $x\le c$,
\cr\vphantom{\vrule width 0pt height 13pt depth 0pt}
-\varepsilon + \frac{h}{\varepsilon}(x-c),&if 
$c\le x \le c+\varepsilon$,
\cr\vphantom{\vrule width 0pt height 13pt depth 0pt}
h-\varepsilon,&if $c+\varepsilon\le x$,
\cr}
$$
let $\tau_{2}(x)=\tau_{1}(x+\varepsilon)+2\varepsilon$,
and let $U$ be the open set defined by (\ref{U}). Then $u$ is a 
Dirichlet $U$-minimizer of $F^{\alpha}$, and a calibration on $U$ is 
given by the piecewise constant function
$$
\phi(x,t):=\cases{(2\lambda, \lambda^{2}),&if 
$c-\varepsilon<x<c+\varepsilon$ and
\cr&
$\varepsilon+\frac{\lambda}2(x-c+\varepsilon)< t < \varepsilon+
\frac{\lambda}2(x-c+\varepsilon)+\frac{\alpha}{2\lambda}$,
\cr\vphantom{\vrule width 0pt height 13pt depth 0pt}
(0,0),&otherwise,
\cr}
$$
where $\lambda>0$ is any constant such that 
$\varepsilon+\varepsilon\lambda +
\frac{\alpha}{2\lambda}\le h-\varepsilon$, for instance 
$\lambda=\sqrt{\alpha}/\sqrt{2\varepsilon}$.
\end{example}

\begin{example}[Harmonic function]\label{ex3} \rom
Let $\Omega$ be a bounded domain in ${\mathbf R}^n$, $n$ arbitrary,
and let $u$ be a harmonic function on $\Omega$. As pointed 
out by Chambolle \cite{Cha}, $u$ is a 
Dirichlet minimizer of $F^{\alpha}$ if
\begin{equation}\label{e1}
\mathop{\rm osc}_{\Omega}u \, 
\mathop{\smash\sup\vphantom\inf}_{\Omega}|\nabla u| \,\le\, \alpha \,,
\end{equation}
where $\mathop{\rm osc}_{\Omega}u:=
\mathop{\smash\sup\vphantom\inf}_{\Omega}u 
-\inf_{\Omega}u$. 
Note that for $n=1$ this condition reduces to the constraint 
$a\lambda^{2}\le \alpha$ of Example~\ref{ex1}. Inspired by the one 
dimensional case (see (\ref{6})), we construct the calibration
\begin{equation}\label{cal}
\phi(x,t):=\cases{
\big( 2\nabla u(x), |\nabla u(x)|^2 \big),&if 
$\frac{1}{2}(u(x)+m)\le t \le \frac{1}{2}(u(x)+M)$,
 \cr\vphantom{\vrule width 0pt height 13pt depth 0pt}
(0,0),&otherwise,
\cr}
\end{equation}
where $m:=\inf_{\Omega}u$ and 
$M:= \mathop{\smash\sup\vphantom\inf}_{\Omega}u$.
Another calibration (see (\ref{6b})) is given by
$$ 
\phi(x,t):=\cases{\big(2\frac{t-m}{u(x)-m}\nabla u(x), 
(\frac{t-m}{u(x)-m})^{2}|\nabla u(x)|^{2}\big),&if
$m\le t\le u(x)$,
\cr\vphantom{\vrule width 0pt height 13pt depth 0pt}
\big(2\frac{M - t}{M- u(x)}\nabla u(x), 
(\frac{M - t}{M- u(x)})^{2}|\nabla u(x)|^{2}\big),&if
$u(x)\le t\le M$,
\cr\vphantom{\vrule width 0pt height 13pt depth 0pt}
(0,0),&otherwise.
\cr}
$$ 
\par
If (\ref{e1}) is not satisfied, $u$ is still
is a Dirichlet $U$-minimizer of $F^{\alpha}$, for
\begin{equation}\label{Unabla}
\textstyle
U:=\{(x,t)\in \Omega\times{\mathbf R}: u(x)-\frac{\alpha}{4}|\nabla u(x)|^{-1}< t
< u(x)+\frac{\alpha}{4}|\nabla u(x)|^{-1} \}\,,
\end{equation}
and a calibration in $U$ is given by $\phi(x,t):=
\big( 2\nabla u(x), |\nabla u(x)|^2 \big)$.
\end{example}

\begin{example}[Pure jump]\label{ex4} \rom
Let $n\ge 2$ and let $\Omega:={]0,a[}\times V$, where $V$ is a 
bounded domain in ${\mathbf R}^{n-1}$ with Lipschitz boundary.
Denoting the first coordinate of $x$ by $x_{1}$,
let $u(x):=0$ for $0<x_{1}<c$, and $u(x):=h$ for
$c<x_{1}<a$, with $0<c<a$ and $h>0$. Using the results of
Example~\ref{ex2} it is easy to see that
$u$ is a Dirichlet minimizer of $F^{\alpha}$ if 
$a\alpha \le h^{2}$. In this case two different calibrations can be 
constructed in the following way:
the projection of
these calibrations onto the $(x_{1},t)$-plane are given by 
(\ref{6}) and (\ref{6b}), with $\lambda=\sqrt\alpha/\sqrt{a}$ and
$x$ replaced by $x_{1}$, while all other components of these 
calibrations vanish.
\par
If $a\alpha > h^{2}$, it may happen that $u$ is still a Dirichlet 
minimizer of $F^{\alpha}$. 
For instance, if $n=2$ and $V={]0,b[}\,$, with $b\alpha\pi  \le 
2h^{2}$, a different calibration has been constucted in 
\cite{ABD}. Therefore $u$ is a Dirichlet minimizer of $F^{\alpha}$ 
even if $a\alpha $ is very large with respect to $h^{2}$, provided 
that $b\alpha$ is small enough.
\par
Arguing as in the last part of Example~\ref{ex2} one can prove that for every $a$ and 
$V$ there exists an open set $U$ of the form (\ref{U}), containing
${\mit\Gamma}_{u}$, such that $u$ 
is a Dirichlet $U$-minimizer of $F^{\alpha}$.
\end{example}

\begin{example}[Triple junction]\label{ex5} \rom
Let $n:=2$, let $\Omega:=B(0,r)$ be the open ball with radius $r>0$ 
centered at the origin, 
and let $u$ be given, in polar coordinates, by
$u(\rho,\theta):=a$  for $0\le\theta<{2\over 3}\pi$, 
$u(\rho,\theta):=b$  for ${2\over 3}\pi\le\theta<{4\over 3}\pi$, 
and $u(\rho,\theta):=c$ for ${4\over 3}\pi\le\theta<2\pi$, 
where $a$, $b$, and $c$ are distinct constants. 
Thus $S_{u}$ is given by three line segments meeting 
at the origin with equal angles.
If
\begin{equation}\label{erre}
2\alpha r \le \min\{|a-b|^{2}, |b-c|^{2}, |c-a|^{2}\}\,,
\end{equation}
then $u$ is a Dirichlet minimizer of $F^{\alpha}$.
To construct a calibration, it is not restrictive to
assume $a<b=0<c$.
Inspired by the one dimensional case 
described in Example~\ref{ex2}, we take 
$e_\pm:=(\pm\sqrt 3/2,-1/2)$,
and $\lambda>0$ such that ${\lambda r \over 2}+
{\alpha \over\lambda}\le \min\{{-a},c\}$ (which is possible by 
(\ref{erre})), 
and we define the calibration by
\begin{equation}\label{tri}
\phi(x,t):=\cases{
(\lambda e_+, \lambda^2/4),&if 
$\frac{\lambda}{4}(r+x\cdot e_+)\le t
 \le \frac{\lambda}{4}(r+x\cdot e_+) + \frac{\alpha}{\lambda}$,
 \cr\vphantom{\vrule width 0pt height 13pt depth 0pt} 
(\lambda e_-, \lambda^2/4),&if
$\frac{\lambda}{4}(-r+x\cdot e_-) - \frac{\alpha}{\lambda} \le t
\le \frac{\lambda}{4}(-r+x\cdot e_-)$,
\cr\vphantom{\vrule width 0pt height 13pt depth 0pt}
(0,0),& otherwise. 
\cr}
\end{equation}
\par
If $\alpha r$ is much larger than $\min\{|a-b|^{2}, |b-c|^{2}, |c-a|^{2}\}$, 
it is easy to construct a comparison function $v$ with the same 
boundary values as $u$ and such that $F^{\alpha}(v)<F^{\alpha}(u)$. 
This shows that in this case $u$ is not a Dirichlet minimizer. 
\par
However, for every value of the parameters $\alpha$, $r$, $a$, $b$, 
$c$, one can construct a suitable neighbourhood $U$ of the graph 
${\mit\Gamma}_u$, of the form (\ref{U}), such that a variant 
of (\ref{tri}) is a calibration in $U$, and therefore $u$ is a 
Dirichlet $U$-minimizer of $F^{\alpha}$. We refer to \cite{ABD} for 
the details.
\end{example}

We consider now the functional $F^{\alpha,\beta}_{g}$, with $\beta>0$.

\begin{example}[Solution of the Neumann problem]
\label{ex6}  \rom
Let $\Omega$ be a bounded open set in ${\mathbf R}^{n}$ with boundary of 
class $C^{1,\varepsilon}$ for some $\varepsilon>0$,
and let $u$ be the 
solution of the Neumann problem
\begin{equation}\label{neumann}
\cases{\Delta u =\beta (u-g)&on $\Omega$,
\cr\vphantom{\vrule width 0pt height 13pt depth 0pt}
\frac{\partial u}{\partial \nu}=0&on
$\partial\Omega$,
\cr}
\end{equation}
with $\beta>0$ and $g\in L^{\infty}(\Omega)$.
Assume that condition (\ref{e1}) of Example~\ref{ex3} is satisfied.
Then $u$ is a minimizer of $F^{\alpha,\beta}_{g}$. If the strict 
inequality holds in (\ref{e1}), then $u$ is the unique 
minimizer. 
 A Neumann calibration $\phi(x,t)$ is given by
$$
\cases{\big(0, \beta|\frac{m}{2}-\frac{u(x)}{2}|^{2}
- \beta|\frac{m}{2}+\frac{u(x)}{2}-g(x)|^{2}\big),&if
$t-\frac{u(x)}{2}< \frac{m}{2}$,
\cr\vphantom{\vrule width 0pt height 13pt depth 0pt}
\big(2\nabla u(x), |\nabla u(x)|^{2} - \beta|t-g(x)|^{2}+ 
\beta|t-u(x)|^{2}\big),&if $\frac{m}{2}\le t-\frac{u(x)}{2} \le
\frac{M}{2}$,
\cr\vphantom{\vrule width 0pt height 13pt depth 0pt}
\big(0, \beta|\frac{M}{2}-\frac{u(x)}{2}|^{2}
- \beta|\frac{M}{2}+\frac{u(x)}{2}-g(x)|^{2}\big),&if
$\frac{M}{2}< t-\frac{u(x)}{2}$,
\cr}
$$
where $m:=\inf_{\Omega}u$ and 
$M:= \mathop{\smash\sup\vphantom\inf}_{\Omega}u$.

If (\ref{e1}) is not satisfied, $u$ is still
is a $U$-minimizer of $F^{\alpha,\beta}_{g}$, where $U$ is the open set
defined by (\ref{Unabla}). A Neumann calibration on $U$ is given by
$$
\phi(x,t):=(2\nabla u(x), |\nabla u(x)|^{2} - \beta|t-g(x)|^{2}+ 
\beta|t-u(x)|^{2}) \,.
$$

The hypothesis that $\partial\Omega$ is of class $C^{1,\varepsilon}$ 
is used only to obtain the boundary condition $\rm(c2)$ of 
Theorem~\ref{3.2}, which, in this case, becomes
\begin{equation}\label{c2}
\lim_{y\to x} \nabla u(y)\cdot \nu(x)=0 \quad
\mbox{for } {\cal H}^{n-1}\mbox{-a.e.\ }x\in\partial\Omega\,.
\end{equation}
It is clear that (\ref{c2}) is still true if for ${\cal H}^{n-1}$-a.e.\ 
$x\in\partial\Omega$ there exists an open neighbourhood $V_{x}$ of $x$ in 
${\mathbf R}^{n}$ such that $V_{x}\cap\partial\Omega$ is a manifold of class 
$C^{1,\varepsilon}$ (see \cite[Theorem~7.5.2]{AFP}). Therefore the result of this example is true 
also when $\Omega$ is polyhedral.
\end{example}

In the next examples we construct a calibration for 
$F^{\alpha,\beta}_{g}$ when the parameter $\beta$ is large enough.

\begin{example}[Smooth $g$ and large $\beta$]\label{ex7}  \rom
Let $\Omega$ be a bounded open set in ${\mathbf R}^{n}$ with smooth boundary,
and let $g\in C^{2}(\overline\Omega)$. There exists a constant 
$\beta_{0}\ge 0$, 
depending on $g$ and $\alpha$, such that for every $\beta>\beta_{0}$ 
the 
solution $u$ of the Neumann problem (\ref{neumann}) 
of Example~\ref{ex6} is the unique minimizer of $F^{\alpha,\beta}_{g}$.
A Neumann calibration is constructed in \cite{ABD}.

This shows that the minimizer of $F^{\alpha,\beta}_{g}$ is smooth, 
provided that $g$ is smooth and $\beta$ is large enough. Therefore the 
solution of the image segmentation problem ($n=2$) based on the 
minimization of $F^{\alpha,\beta}_{g}$ has an empty set of segmentation 
lines if the  ``grey level'' function $g$ is smooth and the parameter $\beta$ in the 
fidelity term $\beta\int_{\Omega}|{u-g}|^{2} dx$  is large.
\end{example}

\begin{example}[Function $g$ with only two values]\label{ex8}
 \rom
Let $\Omega$ be an open set in ${\mathbf R}^{n}$ and let $E$ be
a compact set  contained in $\Omega$ 
with boundary of class $C^{2}$. Let $g(x):=a$ for $x\in E$ and 
$g(x):=b$ for $x\in{\Omega\setminus E}$, with $a\neq b$.
There exists a constant $\beta_{0}\ge0$, 
depending on $g$ and $\alpha$, such that for every 
$\beta > \beta_{0}$ the function $u:=g$ is the unique minimizer of 
$F^{\alpha,\beta}_{g}$.
To construct a calibration, it is not restrictive to assume
$a<b$. We take a $C^{1}$ vectorfield 
$v\colon\Omega\to{\mathbf R}^{n}$ with compact support in $\Omega$ such that 
$|v(x)|\le 1$ for every $x\in\Omega$ and $v(x)$ is the outer unit 
normal to $\partial E$ for every $x\in\partial E$. Then we set
$\phi^{x}(x,t)=\sigma(t)v(x)$,
where $\sigma$ is a fixed positive smooth function with integral 
equal to $\alpha$ and support contained in ${]a,b[}\,$. We see 
that conditions $\rm(b1)$, $\rm(b2)$, and $\rm(c2)$ of 
Theorem~\ref{3.2} are 
satisfied by construction. It remains to choose $\phi^{t}$ so that
$\rm(a1)$, $\rm(a2)$, and $\rm(c1)$ hold. Condition $\rm(a2)$ 
forces us to set $\phi^{t}(x,t)=0$ for $t=g(x)$, while $\rm(c1)$
gives $\partial_{t}\phi^{t}(x,t)=-\sigma(t){\rm div}_{x} v(x)$. These two 
conditions determine $\phi^{t}(x,t)$ at every point $(x,t)$.
It is then easy to see that 
$\rm(a1)$ holds if $\beta$ is large enough.
We refer to \cite{ABD} for the details.
\par
This example shows that, if $g\in SBV(\Omega)$ has only two values, 
and $S_{g}$ is smooth enough, then the minimizer of the Mumford-Shah 
functional
$F^{\alpha,\beta}_{g}$ reconstructs $g$ exactly, when $\beta$ is 
large enough.
\end{example}

Recently the following question has been studied by using the 
calibration method: is 
it true that a function $u$ 
is a (Dirichlet) minimizer of $F^{\alpha,\beta}_{g}$, if it
satisfies the Euler-Lagrange equations and
the domain $\Omega$ is {\it sufficiently small\/}? For the moment 
we have only a partial answer. In \cite{DMM} we have considered the case
where $n:=2$ and $S_{u}$ is a line segment joining two points of the 
boundary of $\Omega$. If $u$ satisfies the Euler-Lagrange equations 
for the ``homogeneous functional''
$F^{\alpha}:=F^{\alpha,0}_{g}$, then for every $x_{0}\in S_{u}$ 
there exists an 
open neighbourhood $\Omega_{0}$ of $x_{0}$, contained in $\Omega$, 
such that $u$ is a Dirichlet minimizer of $F^{\alpha}$ in $\Omega_{0}$.
The minimality is proved by constructing a complicated calibration on 
$\Omega_{0}\times{\mathbf R}$.
\par
This result has been extended in \cite{MM} to the case where $S_{u}$ 
is an analytic curve joining two points of $\partial \Omega$. The 
(more difficult) construction of the calibration
presented in this paper shows that 
one can take the same set $\Omega_{0}$ 
for every $x_{0}\in S_{u}$; in other words, one can take as 
$\Omega_{0}$ a suitable tubular neighbourhood of $S_{u}$. 
Moreover, it is proved in \cite{MM} that
an additional condition on $u$ and 
$S_{u}$ implies that $u$ is a Dirichlet $U$-minimizer for a 
suitable open neighbourhood $U$ of the graph ${\mit\Gamma}_{u}$.
A counterexample (where $S_{u}$ is a line segment joining two points of
 $\partial\Omega$)
shows that this is not always true when $u$ is just a solution  
of the Euler-Lagrange equations with $S_{u}\neq\emptyset$, 
in contrast to the case 
 $S_{u}=\emptyset$ (see Example~\ref {ex6}).

\end{document}